\documentclass[a4paper,12pt]{article}
\usepackage{amsmath,amsthm,amssymb}
\usepackage[matrix,arrow,curve]{xy}

\newtheorem{Lemma}{Lemma}%[section]
\newtheorem{Def}{Definition}%[section]
\newtheorem{prop}{Proposition}%[section]
\newtheorem{theorem}{Theorem}%[section]
\newtheorem*{theorem*}{Theorem}%[section]
\newtheorem{coro}{Corollary}%[section]

\newcommand{\Hom}{{\rm Hom}}
\newcommand{\Id}{{\rm Id}}

\DeclareMathOperator{\Aut}{Aut}
\renewcommand{\Im}{{\rm Im\,}}

\newcommand{\TP}{{\rm TrPic}}

\newcommand{\Pic}{{\rm Pic}}
\newcommand{\Out}{{\rm Out}}
\newcommand{\Inn}{{\rm Inn}}
\newcommand{\add}{{\rm add}}

\newcommand{\End}{{\rm End}}

\newcommand{\td}{\widetilde}

\newcommand{\AAA}{\mathcal{A}}
\newcommand{\BB}{\mathcal{B}}
\newcommand{\CC}{\mathcal{C}}
\newcommand{\DD}{\mathcal{D}}
\newcommand{\TT}{\mathcal{T}}

\newcommand{\KK}{\mathcal{K}}
\newcommand{\NN}{\mathcal{N}}
\newcommand{\RR}{\mathcal{R}}
\newcommand{\PP}{\mathcal{P}}
\newcommand{\QQ}{\mathcal{Q}}
\newcommand{\SSS}{\mathcal{S}}
\newcommand{\II}{\mathcal{I}}
\newcommand{\HH}{\mathcal{H}}
\newcommand{\Kb}{\mathcal{K}^{\rm b}_p}

\renewcommand{\le}{\leqslant}
\renewcommand{\ge}{\geqslant}

\newcommand{\ot}{\otimes}
\newcommand{\kk}{\mathbf{k}}

\newcommand{\floor}[1]{\left\lfloor#1\right\rfloor}

\newenvironment{Proof}[1][Proof]{\begin{trivlist}
\item[\hskip \labelsep {\bfseries #1}]}{\flushright
$\Box$\end{trivlist}}

\numberwithin{equation}{section}

\begin{document}
\title{Derived Picard groups of selfinjective Nakayama algebras.}
\author{Yury Volkov and Alexandra Zvonareva}
\date{}
\maketitle
\begin{center}
    \emph{To the memory of Alexander Ivanov}
\end{center}

\begin{abstract}
In our proceeding paper a generating set of the derived Picard group of a selfinjective Nakayama algebra was constructed combining some previous results for Brauer tree algebras and the technique of orbit categories developed there. In this paper we finish the computation of the derived Picard group of a selfinjective Nakayama algebra.
\end{abstract}

\section{Introduction}
The derived Picard group $\TP(A)$ of an algebra $A$ is the group of
standard autoequivalences of the derived category modulo natural
isomorphisms. This group is an interesting invariant of the derived
equivalence connected to various other invariants such as Hochschild
cohomology \cite{Keller1}, the identity component of the algebraic
group of outer automorphisms \cite{H-ZS} and the group of
autoequivalences of the stable category in the case of a
selfinjective algebra.

The derived Picard group was computed in a few cases. If $A$ is
local or commutative with a connected spectrum, then $\TP(A) \simeq
{\rm Pic}(A) \times \mathbb{Z}$ (see \cite{RouZim}, \cite{Yekut}).
In the work of Miyachi and Yekutieli the derived Picard group was
computed for finite dimensional hereditary algebras \cite{MY}. In
the work of Lenzing and Meltzer \cite{LM} the group of
autoequivalences of the bounded derived category of a canonical
algebra was studied and in the work of Broomhead, Pauksztello and
Ploog this group was described for discrete derived categories
\cite{BPP}.

The derived Picard group of a symmetric Nakayama algebra or an algebra derived equivalent to a symmetric Nakayama algebra was studied in many papers. Rouquier and Zimmermann defined a braid group action on the derived category of a Brauer line algebra with multiplicity $1$ (see \cite{RouZim}). Later, Zimmermann investigated the case of arbitrary multiplicity \cite{Zi}. Schaps and Zakay-Illouz defined an action of a braid group corresponding to $\td A_{n-1}$ on the derived category of a Brauer star algebra with arbitrary multiplicity \cite{ShapZa}. Note that the class of Brauer star algebras with arbitrary multiplicity coincides with the class of symmetric Nakayama algebras. Khovanov and Seidel developed the theory of spherical twists in \cite{ST}, it follows from \cite{ST} and \cite{KS} that the braid group action defined in \cite{RouZim} is faithful. In \cite{Intan} Muchtadi-Alamsyah proved that the braid group action defined in \cite{ShapZa} is faithful in the case of multiplicity $1$. In \cite{Zvonareva} the second named author computed the generating set of $\TP(A)$ for a Brauer star algebra with arbitrary multiplicity. In \cite{VolkZvo} both authors extended this result to the class of selfinjective Nakayama algebras.

Let $\kk$ be an algebraically closed field and $A$ be a selfinjective Nakayama algebra over $\kk$. The aim of this work is to describe the derived Picard group of $A$. The paper is organized as follows: in Sections \ref{dPg} and \ref{sde} we collect some preliminary results on standard derived equivalences and the derived Picard group. In Section \ref{SMfDP} we recall the technique of orbit categories or smash products developed in \cite{VolkZvo}, which is used in Section \ref{fa}.

In Section \ref{mr} we formulate and prove our main result: the description of the derived Picard group for selfinjective Nakayama algebras.

Let $m,n,t>0$ be some integers. We suppose that $n$ and $t$ are
coprime. Let $\QQ_{nm}$ be a cyclic quiver with $nm$ vertices. Let
$\II_{nm,tm}$ be an ideal in the path algebra of $\QQ_{nm}$
generated by all paths of length $tm+1$. We denote $\NN_{nm,
tm}:=\kk\QQ_{nm}/\II_{nm,tm}$.

\begin{theorem*}
$1)$ If $m=1$, then $$\TP(\NN_{nm, tm}) = \TP(\NN_{n, t})\cong
C_n\times\SSS_{\floor{\frac{t+n-1}{n}}}(\kk)\times C_{\infty}.$$

$2)$ If $m>1$ and $t=1$, then $$\TP(\NN_{nm, tm}) = \TP(\NN_{nm,
m})\cong\AAA_{m,n}\times \kk^*.$$

$3)$ If $m>1$ and $t>1$, then $$\TP(\NN_{nm, tm})\cong\big(\BB(\td A_{m-1})\rtimes_{\varphi_{m,n}} C_{nm}\big)\times\SSS_{\floor{\frac{t+n-1}{n}}}(\kk)\times C_{\infty}.$$
\end{theorem*}
The corresponding groups are defined in Section \ref{sg}. In Section
\ref{Pg} we describe the Picard group of a selfinjective Nakayama
algebra, in Section \ref{genset} we describe the generating set of
$\TP(A)$ and prove the first part of the theorem. In Sections
\ref{aR} and \ref{ci} we describe an algebra $\RR$, which is derived
equivalent to $\NN_{nm,m}$. This algebra was studied in
\cite{Efimov}, there a faithful action of a braid group was defined,
providing a generalisation of the results from \cite{ST}, using this
result we prove the second part of the theorem in Section
\ref{case1}. In Section \ref{BAff} we recall some facts about the
affine braid group $\mathcal{B}(\td A_{n-1})$ and it's faithful
action on the derived category in some cases. Also one nice
monomorphism of affine braid groups is described in this section. In
Section \ref{fa} we extend the faithfulness to the general case,
this allows us to prove the third part of the theorem in Section
\ref{case>1}.

\textbf{Acknowledgement:} We would like to thank Sem\"en Podkorytov and Mikhail Basok for discussions on braid groups. The work of the first named author was
partially supported by the Russian Foundation for Basic Research
(RFFI grant 13-01-00902) and by the S$\tilde{\rm a}$o-Paulo Research Foundation (proc. 2014/19521-3), the work of the second named author was
partially supported by the Russian Foundation for Basic Research
(RFFI grant 13-01-00902) and by the Chebyshev Laboratory (Department
of Mathematics and Mechanics, Saint Petersburg State University)
under Russian Federation Government grant 11.G34.31.0026.

\section{Preliminaries}
\subsection{Standard derived equivalences and the derived Picard group}\label{dPg}

Throughout this paper $\kk$ is an algebraically closed field and we
simply write $\ot$ instead of $\ot_{\kk}$. Let $A$ and $B$ be
$\kk$-algebras. All modules are right modules if otherwise is not stated. We
denote by $\CC A$ the category of $A$-complexes. By $\KK A$ we
denote the homotopy category of $\CC A$. By $\Kb A$ we denote the
subcategory of $\KK A$ of bounded complexes of finitely generated
projective modules. By $\DD A$ we denote the derived category of
$\KK A$. The algebras $A$ and $B$
are called derived equivalent if their derived categories are
equivalent as triangulated categories. It is well known that if $A$
and $B$ are derived equivalent then there are a complex of $B^{\rm
op}\ot A$-modules $X$ and a complex of $A^{\rm op}\ot B$-modules $Y$
such that $X\ot_A Y\cong B$ in $\DD(B^{\rm op}\ot B)$ and $Y\ot_B
X\cong A$ in $\DD(A^{\rm op}\ot A)$ \cite{Rickard}. In this case $X$
is called a {\it two-sided tilting complex}. The complex $X$ defines
an equivalence $T_X=-\ot_B X:\DD B\rightarrow \DD A$ (which is
quasi-inverse to $T_Y$). It is well-known that the equivalence $T_X$
induces an equivalence from $\Kb B$ to $\Kb A$. The equivalence
$T_X:\Kb B\rightarrow\Kb A$ is defined modulo natural isomorphism by
the isomorphism class of $X$ in $\DD(A^{\rm op}\ot B)$. Moreover, if
$T_X\cong T_{X'}$, then $X\cong X'$ in $\DD(B^{\rm op}\ot A)$.

\begin{Def}{\rm The equivalences of the form $T_X$ are called {\it standard equivalences}.
The set of isomorphism classes of two-sided tilting $A^{\rm op}\ot A$-complexes forms a
group with respect to the operation $\ot_A$. This group is called
the {\it derived Picard group} of $A$ and we denote it by $\TP(A)$.
As it was mentioned above the group $\TP(A)$ can be considered as a
group of standard autoequivalences of $\Kb A$ modulo natural
isomorphisms. }\end{Def}

\subsection{Standard derived equivalences and tilting functors}\label{sde} It is
not very convenient to work with complexes of bimodules. To
circumvent this problem we need the following definitions.

\begin{Def}{\rm A complex $H\in\Kb A$ is called a {\it tilting complex} for $A$ if the following two conditions are satisfied\\
$1.$ $\Kb A(H,H[i])=0$ for all $i\neq 0$;\\
$2.$ the minimal triangulated subcategory of $\Kb A$ which contains
$H$ and is closed under isomorphisms, direct summands and direct
sums contains $A$.

A {\it tilting functor} from $B$ to $\Kb A$ is a pair $(H,\theta)$
where $H$ is a tilting complex for $A$ and
$\theta:B\rightarrow\End_{\Kb A}(H,H)$ is an isomorphism of
algebras. We say that $(H,\theta)$ is isomorphic to $(H',\theta')$
if there is an isomorphism $\xi\in \Kb A(H,H')$ such that
$\xi\theta(b)=\theta'(b)\xi$ for all $b\in B$. Sometimes we
simply write $\theta$ instead of $(H,\theta)$ if $H$ is clear.}
\end{Def}

Every standard equivalence $T_X:\Kb B\rightarrow \Kb A$ defines a
tilting functor $(H_X,\theta_X)$, since $B$ is an object of $\Kb B$
and $B\cong\End_{\Kb B}(B,B)\cong\End_{\Kb A}(T_X(B),T_X(B))$.
Conversely, if $(H,\theta)$ is a tilting functor, then there is a
two-sided tilting complex $X$ such that
$(H_X,\theta_X)\cong(H,\theta)$ (see \cite[Section 9]{Keller}), but
the construction of $X$ is tedious. We want to construct a functor
$F_{\theta}:\Kb B\rightarrow \Kb A$ in such a way that
$F_{\theta}\cong T_X$. This construction was introduced in
\cite{Rickard2} and is described in \cite[Section 3]{VolkZvo} with a
slight simplification. We recall some details of this construction
here for convenience. Denote by $\PP_B$ the category of finitely
generated projective $B$-modules and by $\add H$ the full
subcategory of $\Kb A$ formed by direct summands of finite direct
sums of copies of $H$. Then $\theta$ induces an equivalence
$S_{\theta}:\PP_B\rightarrow \add H$. The construction of
$S_{\theta}$ is described in \cite{VolkZvo}. Note that the
equivalence $S_{\theta}$ defines the tilting functor $\theta$. So
sometimes we will define the tilting functor $(H,\theta)$ by the
description of the equivalence $S_{\theta}$.

Suppose that we have constructed the functor $S_{\theta}$. Let us
consider $U\in\Kb B$. Then $U=\oplus_{i\in\mathbb{Z}}U_i$, where
$U_i\in\PP_B$. The object $V_i:=S(U_i)\in\add H$ ($i\in\mathbb{Z}$)
is defined and has some grading
$V_i=\oplus_{j\in\mathbb{Z}}V_{i,j}$. The last formula defines a
bigrading on
$V:=S_{\theta}(U)=\oplus_{i\in\mathbb{Z}}V_i=\oplus_{i,j\in\mathbb{Z}}V_{i,j}$.
Note that $V_{i,j}=0$ for large enough and small enough $i$ or $j$.
Let $d_U:U\rightarrow U$ be a differential of $U$. Then
$d_U=\sum_{i\in\mathbb{Z}}d_{U,i}$, where
$d_{U,i}\in\Hom_B(U_i,U_{i+1})$. Let us choose some representative
of $S_{\theta}(d_{U,i})$ in $\CC A(V_i,V_{i+1})$ and denote it by
$d_{1,i}$. Note that $d_{1,i}=\sum\limits_{j\in\mathbb{Z}}d_{1,i,j}$
for some $d_{1,i,j}:V_{i,j}\rightarrow V_{i+1,j}$. We define
$\tau_{1,i,j}:=(-1)^{i+j}d_{1,i,j}$ and
$\tau_1=\sum_{i,j\in\mathbb{Z}}\tau_{1,i,j}$. Also for any
$i\in\mathbb{Z}$ we have a differential
$\tau_{0,i}=\sum\limits_{j\in\mathbb{Z}}\tau_{0,i,j}:V_i\rightarrow
V_i$, where $\tau_{0,i,j}=(d_{V_i})_j:V_{i,j}\rightarrow V_{i,j+1}$.
So we have constructed differentials
$\tau_0=\sum\limits_{i\in\mathbb{Z}}\tau_{0,i},\tau_1:V\rightarrow
V$ of degrees $(0,1)$ and $(1,0)$ respectively. These two
differentials can be extended to a sequence $(\tau_i)_{i\ge 0}$ such
that $\tau_i:V\rightarrow V$ is of degree $(i,1-i)$ and
$\sum\limits_{i=0}^l\tau_i\tau_{l-i}=0$ for any $l\ge 0$. Let us now
define $F_{\theta}(U)\in\Kb A$. As an $A$-module $F_{\theta}(U)=V$.
The $l$-th component of $F_{\theta}(U)$ is defined by the formula
$$
\big(F_{\theta}(U)\big)_l=\bigoplus\limits_{i+j=l} V_{i,j}.
$$
The differential of $F_{\theta}(U)$ is the map induced by the map
$\tau=\sum\limits_{i\ge 0}\tau_i:V\rightarrow V$. Consider
$U,U'\in\Kb B$ and $f\in\Kb B(U,U')$. Suppose that $f$ is
represented by $\sum_{i\in\mathbb{Z}}f_i$ ($f_i:U_i\rightarrow
U_i'$) in $\CC B(U,U')$. Let us introduce $V:=S_{\theta}(U)$,
$V':=S_{\theta}(U')$. The modules $V$ and $V'$ are bigraded as
before. Let $(\tau_i)_{i\ge 0}$ and $(\tau_i')_{i\ge 0}$ be the
corresponding sequences from the construction of $F_{\theta}(U)$ and
$F_{\theta}(U')$ respectively. For any $i\in\mathbb{Z}$ we choose
some representative $\alpha_{0,i}\in\CC A(V_i,V_i')$ of
$S_{\theta}(f_i)$. Then we obtain a morphism
$\alpha_0=\sum_{i\in\mathbb{Z}}\alpha_{0,i}:V\rightarrow V'$ of
degree $(0,0)$. It can be extended to a sequence $(\alpha_i)_{i\ge
0}$ such that $\alpha_i:V\rightarrow V'$ is of degree $(i,-i)$ and
$\sum\limits_{i=0}^l\alpha_i\tau_{l-i}=\sum\limits_{i=0}^l\tau_i'\alpha_{l-i}$
for any $l\ge 0$. We define $F_{\theta}(f)\in\Kb
A(F_{\theta}(U),F_{\theta}(U'))$ as the class of the morphism
$\alpha=\sum\limits_{i\ge 0}\alpha_i:V\rightarrow V'$. It can be
shown (see \cite{Rickard2}, \cite{VolkZvo} for details) that
$F_{\theta}$ is a functor and that $F_{\theta}\cong T_X$ if
$(H_X,\theta_X)\cong(H,\theta)$.

\subsection{Smash product and the derived Picard group}\label{SMfDP}

Let $G$ be a group and let $A$ be a $\kk$-algebra. Let us fix some
map $\Delta:G\rightarrow\Aut A$. We write ${}^g\!a$ instead of
$\Delta(g)(a)$ for $g\in G$, $a\in A$.

\begin{Def}{\rm The {\it smash product} of $A$ and $\kk G$ is the algebra $A\# G$ whose elements are the pairs $a\# g$ ($a\in A$, $g\in G$) with multiplication defined by the
formula
$$
(a\# g)(b\# h)=a{}^g\!b\#gh\,\,(a,b\in A,g,h\in G).
$$}
\end{Def}

Let $U\in\Kb A$, $g\in G$. We define an $A$-complex $U\#g$ in the
following way. The set $(U\#g)_n$ is formed by the elements $x\#g$
($x\in U_n$). The differential $d_{U\#g}$ is defined by the formula
$d_{U\#g}(x\#g)=d_U(x)\#g$ ($x\in U$). The $A$-module structure is
defined by the formula $(x\#g)a=x{}^g\!a\#g$ ($x\in U$, $a\in A$). For
a morphism $f\in\CC A(U,V)$ we define a morphism $f\#g\in\CC
A(U\#g,V\#g)$ by the formula $(f\#g)(x\#g)=f(x)\#g$ ($x\in U$). We
define the $A\#G$-module structure on $U\#G=\oplus_{g\in G}U\#g$ by
the formula $(x\# g)(a\# h)=x{}^g\!a\#gh$ ($a\in U$, $a\in A$, $g,h\in
G$). For a morphism $f\in\CC A(U,V)$ we define a morphism
$f\#G\in\CC(A\#G)(U\#G,V\#G)$ by the formula $(f\#G)(x\#g)=f(x)\#g$
($x\in U$, $g\in G$). Note that $f_1\#G$ is homotopic to $f_2\# G$
and $f_1\#g$ is homotopic to $f_2\# g$ for any $g\in G$ if $f_1$ is
homotopic to $f_2$. So we define $f\#G\in\Kb(A\#G)(U\#G,V\#G)$ and
$f\#g\in\Kb A(U\# g, V\# g)$ ($g\in G$) for $f\in\Kb A(U,V)$. It is
easy to see that we have defined a functor $-\#G:\Kb A\rightarrow
\Kb(A\#G)$. Also for $f\in\Kb A(U,U\#g)$ ($g\in G$) we simply write
$f^k$ instead of $(f\#g^{k-1})\dots(f\# g)f$. Let us define for $g\in G$, $U\in\Kb A$ a
morphism $s_g:(U\#g)\#G\rightarrow U\#G$ by the formula
$$
s_g(x\#g\#h)=x\#gh\,\,(x\in U,h\in G).
$$
Note that $s_g\circ ((f\#g)\#G)=(f\#G)\circ s_g$ for $g\in G$ and $f\in\Kb A(U,V)$.

We recall some results from \cite{VolkZvo} on the connection between
$\TP(A)$ and $\TP(A\#G)$. We will only need the case $G=C_t$ for
some $t\in\mathbb{N}$, where by $C_t$ ($t\in\mathbb{N}\cup\{\infty\}$) we denote
the cyclic group of order $t$. We fix some generating element of
$C_t$ for each $t$ and denote it by $r$.

A {\it tilting $C_t$-functor} from $A$ to itself is a triple $(H,\theta,\psi)$, where $(H,\theta)$ is a tilting functor and $\psi\in\Kb A(H, H\# r)$ is an isomorphism which satisfies the following conditions:\\
\begin{itemize}
\item $\psi\theta({}^{r^{-1}}\!\!a)=(\theta(a)\# r)\psi$ for any $a\in A$;
\item $\psi^t=\Id_H$  in $\Kb A$.
\end{itemize}
Note that $\psi^{k+t}=\psi^k$ for any $k\ge 0$ by the last condition
and so we can define $\psi^k$ for any $k\in\mathbb{Z}$. Let now
$(H,\theta,\psi)$ be a tilting $C_t$-functor from $A$ to $A$. We
define an isomorphism $\theta\#\psi:A\# C_t\rightarrow\Kb(A\#
C_t)(H\#C_t,H\#C_t)$ by the formula
$$
(\theta\#\psi)(a\# r^k)=s_{r^{k}}\circ(\psi^{k}\# C_t)\circ(\theta({}^{r^{-k}}\!\!a)\#C_t)\,\,(a\in A, k\in\mathbb{Z}).
$$
Thus, there is a tilting functor $(H\# C_t,\theta\#\psi)$ from $A\# C_t$ to itself.

The following lemma follows from \cite{VolkZvo}.

\begin{Lemma}\label{group_act}
There is a group $\TP_{C_t}(A)$ and two homomorphisms $\Phi_A:\TP_{C_t}(A)\rightarrow \TP(A)$ and $\Psi_A:\TP_{C_t}(A)\rightarrow \TP(A\#C_t)$ which satisfy the following conditions
\begin{itemize}
\item for any tilting $C_t$-functor $(H,\theta,\psi)$ there is an element $W\in\TP_{C_t}(A)$ such that $\Phi_A(W)=F_{\theta}$ and $\Psi_A(W)=F_{\theta\#\psi}$;
\item for any $W\in\TP_{C_t}(A)$ the condition $\Phi_A(W)\in\Pic(A)$ is equivalent to the condition $\Psi_A(W)\in\Pic(A\#C_t)$.
\end{itemize}
\end{Lemma}

\subsection{Some groups}\label{sg}

In this section we define some groups which we will use to formulate
our main result.

The first group we need is the group of automorphisms of the algebra
$\kk[x]/x^{N+1}$ ($N\ge 1$). We denote this group by $\SSS_N(\kk)$.
In more detail, $\SSS_N(\kk)$ is a group whose elements are the
sequences $(a_i)_{1\le i\le N}$ ($a_i\in\kk$, $a_1\neq 0$) and the
group operation is defined by the formula
\begin{equation} \label{mult}
(a_i)_{1\le i\le n}(b_i)_{1\le i\le
n}=\left(\sum\limits_{j=1}^ia_j\sum\limits_{k_1+\dots+k_j=i}b_{k_1}\dots
b_{k_j}\right)_{1\le i\le n}.
\end{equation}
In particular, $\SSS_1(\kk)=\kk^*$ is the multiplicative group of
the field $\kk$.

We recall now the definition of braid groups associated to $A_n$
and $\td A_{n-1}$. These groups are particular cases of a braid
group associated to a Coxeter group or a graph but we are going to
need only these cases.

\begin{Def} {\rm  The  {\it braid group associated to $A_n$} is the group $\BB(A_n)$ defined by generators $s_1,\dots, s_n$ and relations:
$$s_is_{i+1}s_i=s_{i+1}s_is_{i+1} \mbox{ for $1\le i\le n-1$ and } s_is_j=s_js_i \mbox{ for $1\le i,j\le n$, } |i-j|>1.$$}
\end{Def}
\begin{Def} {\rm  The  {\it braid group associated to $\td A_{n-1}$} for $n> 2$ is the group $\BB(\td A_{n-1})$ defined by generators $\{s_i\}_{i\in \mathbb{Z}_{n}}$ and relations:
$$s_is_{i+1}s_i=s_{i+1}s_is_{i+1} \mbox{ for $i\in \mathbb{Z}_{n}$ and } s_is_j=s_js_i \mbox{ for $i,j\in \mathbb{Z}_{n}$, } i-j\neq 1, -1.$$
The  {\it braid group associated to $\td A_{n-1}$} for $n = 2$ is
the group $\BB(\td A_{1})= F_2,$ i.e. the free group with two
generators.}
\end{Def}

Now we define the group $\AAA_{m,n}$ ($m\ge 2$, $n\ge 1$) as the group generated by the elements $s_1,\dots, s_m, r_1, r_2$ satisfying the relations:
$$s_is_{i+1}s_i=s_{i+1}s_is_{i+1} \mbox{ for $1\le i\le m-1$ and } s_is_j=s_js_i \mbox{ for $1\le i,j\le m$, } |i-j|>1,$$ $$s_ir_k=r_ks_i\mbox{ for $1\le i\le m$, $k\in\{1,2\}$}, r_1r_2=r_2r_1, r_2^n=1 \mbox{ and } (s_m\dots s_1)^{m+1}=r_2^{m+1}r_1^{2m}.$$

Let $G,H$ be groups and let $\varphi:H\rightarrow\Aut(G)$ be a
homomorphism from $H$ to the group of automorphisms of $G$. We
denote by $G\rtimes_{\varphi}H$ the semidirect product of $G$ and
$H$ such that
$$(g_1,h_1)(g_2,h_2)=(g_1\varphi(h_1)(g_2),h_1h_2)$$
for $g_1,g_2\in G$, $h_1,h_2\in H$. Let us define the homomorphism
$\varphi_{m,n}:C_{nm}\rightarrow\BB(\td A_{m-1})$ by the formula
$$
\varphi_{m,n}(r^l)(s_i)=s_{i+l}\,\,\,(i\in\mathbb{Z}_m,
l\in\mathbb{Z}).
$$
Thus, the group $\BB(\td A_{m-1})\rtimes_{\varphi_{m,n}}C_{nm}$ is
defined.

\section{Main result}\label{mr}

Let us consider the algebra $\NN_{nm, tm}$ defined in the following
way. Let $m,n,t>0$ be some integers. We suppose that $n$ and $t$ are
coprime. Let $\QQ_{nm}$ be a cyclic quiver with $nm$ vertices, i.e.
the quiver whose vertex set is $\mathbb{Z}_{nm}$ and whose arrows
are $\beta_i:i\rightarrow i+1$ ($i\in\mathbb{Z}_{nm}$). Let
$\II_{nm,tm}$ be an ideal in the path algebra of $\QQ_{nm}$
generated by all paths of length $tm+1$. We denote $\NN_{nm,
tm}:=\kk\QQ_{nm}/\II_{nm,tm}$.

 $$
 \QQ_{nm}:
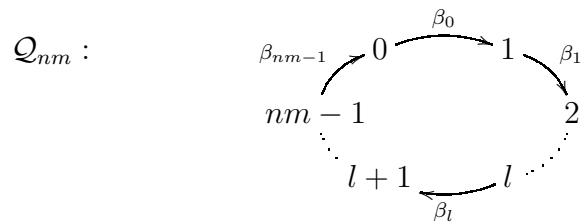
\hspace{2cm}
\xymatrix@!0
{
&0\ar@/^/[rr]^{\beta_0}&&1\ar@/^/[rd]^{\beta_1}\\
nm-1\ar@/^/[ru]^{\beta_{nm-1}}&&&&2\ar@/^/@{.}[ld]\\
&l+1\ar@/^/@{.}[lu]&&l\ar@/^/[ll]^{\beta_l}
}
$$

In this work we prove the following theorem.

\begin{theorem}\label{MainThm}
$1)$ If $m=1$, then $$\TP(\NN_{nm, tm}) = \TP(\NN_{n, t})\cong
C_n\times\SSS_{\floor{\frac{t+n-1}{n}}}(\kk)\times C_{\infty}.$$

$2)$ If $m>1$ and $t=1$, then $$\TP(\NN_{nm, tm}) = \TP(\NN_{nm,
m})\cong\AAA_{m,n}\times \kk^*.$$

$3)$ If $m>1$ and $t>1$, then $$\TP(\NN_{nm, tm})\cong\big(\BB(\td A_{m-1})\rtimes_{\varphi_{m,n}} C_{nm}\big)\times\SSS_{\floor{\frac{t+n-1}{n}}}(\kk)\times C_{\infty}.$$
\end{theorem}

The proof of the first part of the theorem can be found in Section
\ref{genset}, the proof of the second part of the theorem can be
found in Section \ref{case1}, the proof of the third part of the
theorem can be found in Section \ref{case>1}.

From here on we fix some integers $m,n,t>0$. We omit them in some
notation which in fact depends on them. For example, we simply write
$\NN$ instead of $\NN_{nm,tm}$.

\subsection{The description of the Picard group}\label{Pg}

Denote by $\Pic(A)$ the Picard group of an algebra $A$, i.e. the
group of isomorphism classes of invertible $A^{op}\otimes A$-modules
or equivalently the group of Morita autoequivalences of $A$ modulo
natural isomorphisms. If $(H,\theta)$ is a tilting functor from $B$
to $\Kb A$, the equivalence $F_{\theta}$ is a Morita equivalence iff
$H$ is isomorphic to some object $H'\in \Kb A$ concentrated in
degree zero ($H'$ is concentrated in degree zero if $H'_i=0$ for
$i\neq 0$). The group ${\rm Out}(A) = {\rm Aut}(A) /{\rm Inn}(A)$ of
outer automorphisms of $A$ coincides with ${\rm Pic}(A)$
\cite{Bolla} and is clearly a subgroup of $\TP(A)$. Every
automorphism $\sigma$ of $A$ determines a tilting functor
$(A,\bar\sigma)$, where $\bar\sigma(a)(a')=\sigma(a)a'$ for $a,a'\in
A$. A tilting functor $(A,\bar\sigma)$ is isomorphic to
$(A,\bar\sigma')$ if and only if $\sigma$ coincides with $\sigma'$
modulo the subgroup of inner automorphisms. For $\sigma\in\Aut(A)$
we write simply $\sigma$ instead of $F_{\bar\sigma}$. Also we denote
by $\Pic_0(A)$ the subgroup of $\TP(A)$ generated by standard
equivalences $F_{\theta}$ such that $F_{\theta}(P)\cong P$ in $\Kb
A$ for any $P\in\PP_A$. Note that $\Pic_0(A)\cong\Pic_0(\NN)$ if $A$
is derived equivalent to $\NN$, since $\NN$ is a uniserial algebra
(see \cite[Section 4]{Linck}).
Note that if $X$, $X'$ are two two-sided tilting complex of $B^{\rm
op}\ot A$-modules such that $X\cong X'$ in $\DD(A)$, then there exists $\sigma \in \Pic(B)$ such that $X\cong \sigma X'$ in $\DD(B^{\rm op}\ot A)$ \cite[Proposition 2.3]{RouZim}. If the action of the equivalences induced by $X$, $X'$ agree on projective modules, then $\sigma \in \Pic_0(B)$.

For $i\in\mathbb{Z}_{nm}$ we denote by $e_i$ the primitive
idempotent of $\NN$ corresponding to the vertex $i$ and by $P_i$ the
projective module $e_i\NN$. For $w\in e_j\NN e_i$ we denote by $w$
the unique homomorphism from $P_i$ to $P_j$ which sends $e_i$ to
$w$. Also we introduce the following auxiliary notation:
$$
\beta_{i,k}=\beta_{i+k-1}\dots\beta_i.
$$
In particular, $\beta_{i,0}=e_i$. We write simply $\beta$ instead of $\beta_i$ and $\beta^k$ instead of $\beta_{i,k}$ if $i$ is clear.

\begin{Def}{\rm The {\it rotation of $\NN$} is an automorphism $\rho\in\Aut(\NN)$ defined by the formulas
$\rho(e_i)=e_{i+1}$ and $\rho(\beta_i)=\beta_{i+1}$.
}
\end{Def}

Let us consider a sequence $c=(c_i)_{1\le i\le N}$, where
$c_1,\dots,c_N\in\kk$ and $c_1\neq 0$. Let us introduce the notation
$N:=\floor{\frac{t+n-1}{n}}$,
$u_c:=\sum\limits_{i=1}^Nc_i\beta_{1,(i-1)nm}$. We denote by
$\mu_c\in\Aut(\NN)$ the automorphism of $\NN$ defined by the
formulas
$$\mu_c(e_i)=e_i\,\,(i\in\mathbb{Z}_{nm}),\,\,\mu_c(\beta_i)=\beta_{i}\,\,(i\in\mathbb{Z}_{nm}\setminus\{0\})\mbox{ and }\mu_c(\beta_0)=u_c\beta_0.$$

Let us define the map $\xi:C_{nm}\times\SSS_N(\kk)\rightarrow {\rm
Out}(\NN)$ by the formula
 $
 \xi(r^l,c)= [\rho^l\mu_c],
 $
where $[f]$ denotes the class of an automorphism $f$ in ${\rm Out}(\NN)$.
 \begin{prop}\label{Picard} $\xi$ is an isomorphism of groups.
 \end{prop}
\begin{Proof} Let us prove that $\xi$ is a homomorphism. It is not hard to check that $\mu_c\mu_{c'}=\mu_{cc'}$, where the product $cc'$ is defined by the formula \eqref{mult}. Since $\rho^{nm}=\Id_{\NN}$, it remains to prove that $[\rho\mu_c]=[\mu_c\rho]$ for any $c\in\SSS_N(\kk)$. Let us consider the automorphism $\sigma:\NN\rightarrow\NN$ defined by the formula $\sigma(x)=(1-e_1+u_c)^{-1}x(1-e_1+u_c)$ for $x\in\NN$. It is easy to check that $\rho\mu_c=\sigma\mu_c\rho$. Since the automorphism $\sigma$ is inner, we have proved that $\xi$ is a homomorphism.

Let us now prove that $\xi$ is injective. We have to prove that if
$[\rho^l\mu_c]=[\Id_{\NN}]$, then $nm\mid l$ and
$c=1_{\SSS_N(\kk)}$. Assume that there is $y\in\NN^*$ such that
$\rho^l\mu_c(x)=y^{-1}xy$ for any $x\in\NN$. The element $y \in \NN$
is of the form
$$y = \sum_{i\in\mathbb{Z}_{nm}} (d_{i,0}e_i +
d_{i,1}\beta_{i,1} + d_{i,2}\beta_{i,2} + \dots +
d_{i,tm}\beta_{i,tm}), $$ where $d_{i,0} \in \kk^*, d_{i,k} \in \kk$
($i\in\mathbb{Z}_{nm}$, $1 \leq k \leq tm$). Since
$y\rho^l\mu_c(e_0)=ye_l=d_{l,0}e_l + d_{l,1}\beta_{l,1}+\dots$ and
$e_0y=d_{0,0}e_0 + d_{nm-1,1}\beta_{nm-1,1}+\dots$, we have $nm\mid l$ and so $\rho^l=\Id_{\NN}$.

We deduce from the equalities $y\beta_i=\beta_iy$
($i\in\mathbb{Z}_{nm}\setminus\{0\}$) that $d_{i+1,knm}=d_{i,knm}$
for $i\in\mathbb{Z}_{nm}\setminus\{0\}$, $0\le k\le N-1$.
Consequently, $d_{1,knm}=d_{2,knm}=\dots=d_{0,knm}$ for $0\le k\le
N-1$. We deduce from the last equality and the equality
$yu_c\beta_0=\beta_0y$ that
\begin{equation}\label{eq_35}
\sum_{i=0}^{k}d_{0,inm}c_{k-i+1}=d_{0,knm}\,\,(0\le k\le N-1).
\end{equation}
Since $d_{0,0}\neq 0$, it follows from \eqref{eq_35}
for $k=0$ that $c_1=1$. Suppose that we have proved that $c_1=1$ and
$c_i=0$ for $2\le i\le s$ ($s<N$). Then it follows from
\eqref{eq_35} for $k=s$ that $c_{s+1}=0$. So we have  $c_1=1$ and
$c_i=0$ for $2\le i\le N$, i.e. $c=1_{\SSS_N(\kk)}$. It remains to
prove that $\xi$ is surjective.

If $A$ is a finite-dimensional split algebra over a field, it admits
a Wedderbern-Malcev decomposition, i.e. $A=B \oplus J$, where $J$ is
the Jacobson radical of $A$, let us fix such a decomposition.
Following \cite{Pollack} and \cite{Guil-AsencioSaorin} we can
identify $\Out(A)$ with a group $\widehat{H}_A/\widehat{H}_A \cap
\Inn(A)$, where $\widehat{H}_A = \{f \in \Aut(A) \mid f(B) \subseteq
B\}$. For $A=\NN$ fix $B=\langle e_o, \dots
e_{nm-1}\rangle_\kk$. Thus any outer automorphism has a
representative $f$ such that $f(B) \subseteq B$.

Let us prove that the class in $\Out(\NN)$ of any automorphism $f
\in \widehat{H}_{\NN}$ lies in the image of $\xi$. It is easy to see
that $f$ is given by the following formulas:
$$f(e_i)=e_{i+l}, \mbox{ } f(\beta_i)=b_{i,1}\beta_{i+l,1} + b_{i,2}\beta_{i+l,nm+1} + \dots + b_{i,N}\beta_{i+l,(N-1)nm+1},$$
where $0\le l< nm$, $b_{i,1}\in\kk^*$, $b_{i,k}\in\kk$ for
$i\in\mathbb{Z}_{nm}$, $2\le k\le N$. Applying
$\rho^{-l}=\xi(r^{-l})$ we can assume that $l=0$. Let us construct
by induction an automorphism $g_k\in\Im\xi$ ($0\le k\le nm$) such that
$g_kf(e_i)=e_i$ ($i\in\mathbb{Z}_{nm}$) and $g_kf(\beta_i)=\beta_i$
($0\le i< k$). It is clear that we can take $g_0=\Id_{\NN}$. Suppose
that we have constructed the required automorphism $g_k$ for some
$0\le k< nm$. We have $$g_kf(\beta_k)=a_{k,1}\beta_{k,1} +
a_{k,2}\beta_{k,nm+1} + \dots +
a_{k,N}\beta_{i,(N-1)nm+1}=\rho^k(u_{a_k})\beta_k$$ for some
$a_k=(a_{k,i})_{1\le i\le N}\in\SSS_N(\kk)$. In this case we set
$g_{k+1}=(\rho^k\mu_{a_k}^{-1}\rho^{-k})g_k$. It is easy to check
that $g_{k+1}$ satisfies the required conditions. Then
$g_{nm}f=\Id_{\NN}$, i.e. $f=g_{nm}^{-1}\in\Im\xi$.
 \end{Proof}

 \begin{coro}\label{Pic0}
$ \Pic_0(\NN)\cong\SSS_N(\kk).$
 \end{coro}

\subsection{Generating set for $\TP(\NN)$}\label{genset}

Let $\NN$ be as above. If $U$ is a module, then we also denote by
$U$ the corresponding complex concentrated in degree 0. When we
write $\dots\rightarrow\underline{U} \rightarrow \dots$ we mean the
complex whose zero component is $U$. If $V,V'$ are complexes, then
by $(\dots,\underline{f},\dots)$ we denote the morphism of complexes
whose zero component is $f:V_0\rightarrow V_0'$. We simply write $f$
instead of $(\dots,0,\underline{f},0,\dots)$. For $i\in\mathbb{Z}_{nm}$, $1\le k\le
m-1$ we introduce the complexes
$$
X_i:=P_{i-tm}\xrightarrow{\beta} P_{i-tm+1}\xrightarrow{\beta^{tm}} \underline{P_{i+1}}\mbox{ and }Y_{i,k}:=\underline{{P_i}}\xrightarrow{\beta^k} P_{i+k}.
$$

If $m>1$, then for $l\in\mathbb{Z}_m$ we introduce the $\NN$-complex
$$
H_l=\Big(\bigoplus_{i\in\mathbb{Z}_{nm}, m\nmid i-l}P_i\Big)\oplus\Big(\bigoplus_{i\in\mathbb{Z}_{nm}, m\mid i-l}X_i\Big).
$$
In this case we define an equivalence
$S_{\theta_l}:\PP_{\NN}\rightarrow \add H_l$ in the following way:
$$
S_{\theta_l}(P_i)=\begin{cases}
P_i,&\mbox{ if $m\nmid i-l$ and $m\nmid i-1-l$,}\\
P_{i+1},&\mbox{ if $m\mid i-l$,}\\
X_{i-1},&\mbox{ if $m\mid i-1-l$,}
\end{cases}
$$
$$
S_{\theta_l}(\beta_i)=\begin{cases}
\beta_i,&\mbox{ if $m\nmid i-l+1$ and $m\nmid i-1$,}\\
\beta_{i+1}\beta_i,&\mbox{ if $m\mid i-l+1$,}\\
\Id_{P_{i+1}},&\mbox{ if $m\mid i-l$.}
\end{cases}
$$

If $m>1$ and $t=1$, then for $l\in\mathbb{Z}_m$ we introduce the $\NN$-complex
$$
Q_l=\bigoplus_{i\in\mathbb{Z}_{nm}, m\mid i-l}\Big(P_i\oplus \bigoplus_{k=1}^{m-1}Y_{i,k}\Big).
$$
In this case we define an equivalence $S_{\varepsilon_l}:\PP_{\NN}\rightarrow \add Q_l$ in the following way:
$$
S_{\varepsilon_l}(P_i)=\begin{cases}
P_{i-m},&\mbox{ if $m\mid i-l$,}\\
Y_{i+k-m,m-k},&\mbox{ if $m\mid i+k-l$ for some $k$, $1\le k\le
m-1$,}
\end{cases}
$$
$$
S_{\varepsilon_l}(\beta_i)=\begin{cases}
\beta_{i-m,m},&\mbox{ if $m\mid i-l$,}\\
\Id_{P_{i-m+1}},&\mbox{ if $m\mid i-l+1$,}\\
(\underline{\Id_{P_{i+k-m}}},\beta_i),&\mbox{ if $m\mid i-l+k$ for some $k$, $2\le k\le m-1$.}
\end{cases}
$$
For simplicity of notation we write $\HH_l$ instead of
$F_{\theta_l}$ and $\QQ_l$ instead of $F_{\varepsilon_l}$.

The following theorem was proved in \cite{VolkZvo}.

\begin{theorem}\label{Nakayama_gen}
$1)$ If $m=1$, then $\TP(\NN)$ is generated by the shift and $\Pic(\NN)$;\\
$2)$ If $m>1$, $t>1$, then $\TP(\NN)$ is generated by the shift, $\Pic(\NN)$ and $\HH_l$ ($l\in\mathbb{Z}_m$);\\
$3)$ If $m>1$, $t=1$, then $\TP(\NN)$ is generated by the shift, $\Pic(\NN)$, $\HH_l$ and $\QQ_l$ ($l\in\mathbb{Z}_m$).
\end{theorem}

\begin{Proof}[Proof of part 1 of Theorem \ref{MainThm}] Let $m=1$. By part 1 of Theorem \ref{Nakayama_gen} we have $\TP(\NN)=\langle \Pic(\NN),G\rangle$, where $G\cong C_{\infty}$ is the cyclic group generated by the shift. By Proposition \ref{Picard} we have $\Pic(\NN)\cong C_n\times\SSS_{\floor{\frac{t+n-1}{n}}}(\kk)$.
It is clear that $\Pic(\NN)\cap G=\{\Id_{\Kb\NN}\}$. Since the shift commutes with any element of $\Pic(\NN)$, we have
$$\TP(\NN)=\Pic(\NN)\times G\cong C_n\times\SSS_{\floor{\frac{t+n-1}{n}}}(\kk)\times C_{\infty}.$$
\end{Proof}

From here on we assume that $m\ge 2$. We need the following result (which is similar to \cite[Proposition 3]{ShapZa}).

\begin{Lemma}\label{aff_braid}
There is a homomorphism $\eta:\BB(\td A_{m-1})\rightarrow \TP(\NN)$
which sends $s_i$ to $\HH_i$ for all $i\in\mathbb{Z}_m$.
\end{Lemma}
\begin{Proof} If $m=2$, then $\BB(\td A_{m-1})=F_2$ and there is nothing to check. Assume that $m>2$. We need to check that $\HH_i$ satisfy the relations of
the braid group associated to $\td A_{m-1}$. The relations
$\HH_i\HH_j=\HH_j\HH_i \mbox{ for } |i-j|>1$ are straightforward.
Let us check that $\HH_l\HH_{l+1}\HH_l=\HH_{l+1}\HH_l\HH_{l+1}.$ For
simplicity let us assume that $m>3$, the case $m=3$ is similar. The diagram below shows how $\HH_l$, $\HH_{l+1}\HH_l$ and
$\HH_l\HH_{l+1}\HH_l$ act on the projective modules and morphisms
between them. In this computation we can take $\tau_i=0,$ $i>1$. The
morphism $\HH_{l+1}\HH_l(\beta_{km+l+1})$ can be computed using the
construction from Section \ref{sde} with $\alpha_1=-\Id,$
$\alpha_2=-\beta_{km+l}$, in all other cases $\alpha_i=0,$ $i>0$.
The action of $\HH_l\HH_{l+1}\HH_l$ on the indecomposable projective
modules and morphisms between them not appearing in this diagram is
trivial. {\scriptsize $$ \xymatrix @R=2pc @C=1.5pc {
P_{km+l-1} \ar[d]^{\beta}&&&P_{km+l-1}\ar[d]^{\beta^2}&&& P_{km+l-1}\ar[d]^{\beta^3}\\
P_{km+l}\ar[d]^{\beta}&&&P_{km+l+1}\ar[d]^{\Id}&&& P_{km+l+2} \ar[d]^{\Id} \\
P_{km+l+1} \ar[d]^{\beta}  \ar@{~>}[r]^{\HH_l} &P_{(k-t)m+l}\ar[r]^{\beta}&P_{(k-t)m+l+1}\ar[r]^{\beta^{tm}}&P_{km+l+1}\ar@{~>}[r]^{\HH_{l+1}}  \ar[d]^{\beta}
&P_{(k-t)m+l} \ar[r]^{\beta^{2}} \ar[d]^{-\beta}&P_{(k-t)m+l+2} \ar[r]^{-\beta^{tm}} \ar[d]^{-\Id} & P_{km+l+2} \ar[d]^{\Id} \\
P_{km+l+2}\ar[d]^{\beta}&&&P_{km+l+2} \ar[d]^{\beta} &P_{(k-t)m+l+1} \ar[r]^{\beta}&P_{(k-t)m+l+2} \ar[r]^{\beta^{tm}} & P_{km+l+2} \ar[d]^{\beta} \\
P_{km+l+3}&&&P_{km+l+3}&&& P_{km+l+3}\\
}
$$}
{\scriptsize $$
\xymatrix @R=2pc @C=1.5pc {
&&&&&P_{km+l-1} \ar[d]^{\beta^3}\\
&&&&&P_{km+l+2} \ar[d]^{\Id}\\
\ar@{~>}[r]^{\HH_l}&&&P_{(k-t)m+l+1} \ar[r]^{\beta} \ar[d]^{-\Id}&P_{(k-t)m+l+2} \ar[r]^{\beta^{tm}} \ar[d]^{-\Id} & P_{km+l+2} \ar[d]^{\Id}\\
&P_{(k-2t)m+l} \ar[r]^{\beta}&P_{(k-2t)m+l+1} \ar[r]^{\beta^{tm}}&P_{(k-t)m+l+1} \ar[r]^{\beta}&P_{(k-t)m+l+2} \ar[r]^{-\beta^{tm}} & P_{km+l+2} \ar[d]^{\beta}\\
&&&&&P_{km+l+3}\\
}
$$}
Let us compute $\HH_l\HH_{l+1}.$
{\scriptsize $$
\xymatrix @R=2pc @C=1.5pc {
P_{km+l-1} \ar[d]^{\beta}&&&&& P_{km+l-1}\ar[d]^{\beta^2}&\\
P_{km+l}\ar[d]^{\beta}&&&&&P_{km+l+1}\ar[d]^{\beta}&\\
P_{km+l+1} \ar[d]^{\beta}  \ar@{~>}[r]^(.6){\HH_l\HH_{l+1}}&&&&&P_{km+l+2}\ar[d]^{\Id}\ar@{~>}[r]^(.7){\HH_{l+1}}&\\
P_{km+l+2}\ar[d]^{\beta}&P_{(k-2t)m+l}\ar[r]^{\beta}&P_{(k-2t)m+l+1}\ar[r]^{\beta^{tm}}&P_{(k-t)m+l+1}\ar[r]^{\beta}&P_{(k-t)m+l+2}\ar[r]^{-\beta^{tm}}&P_{km+l+2}\ar[d]^{\beta}&\\
P_{km+l+3}&&&&&P_{km+l+3}&\\
}
$$}
$\HH_{l+1}\HH_l\HH_{l+1}(P_{km+l+2})$ is the totalization of the following bicomplex, where $\tau_2$ and $\tau_3$ are constructed according to the algorithm from Section \ref{sde} and are shown as dotted lines.
{\scriptsize $$
\xymatrix @R=2pc @C=1.5pc {
P_{(k-2t)m+l}\ar[r]^{\beta^2} \ar@{-->}[ddrrr]^{-\beta}&P_{(k-2t)m+l+2} \ar@{-->}[ddrrr]_(.4){-\beta^{tm-1}} \ar[r]^{-\beta^{tm}} \ar@{-->}[drr]^(.7){\Id}&P_{(k-t)m+l+2}\ar[r]^{\Id} \ar@{-->}[drr]^(.7){-\Id}&P_{(k-t)m+l+2}\ar[r]^{\beta^{tm}}&P_{km+l+2}\\
&&&P_{(k-2t)m+l+2}\ar[u]^(.3){\beta^{tm}}&P_{(k-t)m+l+2}\ar[u]^(.3){\beta^{tm}}\\
&&&P_{(k-2t)m+l+1}\ar[u]^{\beta}&P_{(k-t)m+l+1}\ar[u]^{\beta}\\
}
$$}
{\footnotesize
$$\HH_{l+1}\HH_l\HH_{l+1}(P_{km+l+2})=P_{(k-2t)m+l}\xrightarrow{\left(\begin{smallmatrix}
\beta^2\\-\beta \end{smallmatrix}\right)} P_{(k-2t)m+l+2} \oplus
P_{(k-2t)m+l+1}\xrightarrow{\left(\begin{smallmatrix} -\beta^{tm} &
0\\\Id &\beta\\-\beta^{tm-1}&0 \end{smallmatrix}\right)}$$
$$P_{(k-t)m+l+2} \oplus P_{(k-2t)m+l+2} \oplus
P_{(k-t)m+l+1}\xrightarrow{\left(\begin{smallmatrix} \Id& \beta^{tm}
& 0\\-\Id &0&\beta \end{smallmatrix}\right)} P_{(k-t)m+l+2} \oplus
P_{(k-t)m+l+2}\xrightarrow{\left(\begin{smallmatrix} \beta^{tm} &
\beta^{tm} \end{smallmatrix}\right)} P_{km+l+2}$$}

The isomorphism
from $\HH_{l+1}\HH_l\HH_{l+1}(P_{km+l+2})$ to
$\HH_{l}\HH_{l+1}\HH_{l}(P_{km+l+2})$ is given by
$$f=(\Id,
\left(\begin{smallmatrix} 0 & -\Id \end{smallmatrix}\right),
\left(\begin{smallmatrix} 0 & -\beta^{tm-1}&-\Id
\end{smallmatrix}\right), \left(\begin{smallmatrix} -\Id & -\Id
\end{smallmatrix}\right), \underline{\Id})$$
with the inverse
$$f^{-1}=\left(\Id,
\left(\begin{smallmatrix} \beta \\ -\Id \end{smallmatrix}\right),
\left(\begin{smallmatrix} 0 \\ 0\\-\Id \end{smallmatrix}\right),
\left(\begin{smallmatrix} 0 \\ -\Id \end{smallmatrix}\right),
\underline{\Id}\right).$$
  Indeed, it is easy to see that $\HH_{l}\HH_{l+1}\HH_{l}(P_{km+l+2})$ is a direct summand of $\HH_{l+1}\HH_l\HH_{l+1}(P_{km+l+2})$,
since $\HH_{l+1}\HH_l\HH_{l+1}(P_{km+l+2})$ is an image of an
indecomposable object under an equivalence, the other direct summand
is isomorphic to zero. And we get: {\scriptsize $$ \xymatrix @R=2pc
@C=1.5pc {
&&&&&P_{km+l-1} \ar[d]^{\beta^3}\\
&&&&&P_{km+l+2} \ar[d]^{\Id}\\
\ar@{~>}[r]^{\HH_{l+1}}&&&P_{(k-t)m+l+1} \ar[r]^{\beta} \ar[d]^{-\Id}&P_{(k-t)m+l+2} \ar[r]^{\beta^{tm}} \ar[d]^{-\Id} & P_{km+l+2} \ar[d]^{\Id}\\
&P_{(k-2t)m+l} \ar[r]^{\beta}&P_{(k-2t)m+l+1} \ar[r]^{\beta^{tm}}&P_{(k-t)m+l+1} \ar[r]^{\beta}&P_{(k-t)m+l+2} \ar[r]^{-\beta^{tm}} & P_{km+l+2} \ar[d]^{\beta}\\
&&&&&P_{km+l+3}\\
}
$$}
Here we replace $\HH_{l+1}\HH_l\HH_{l+1}(\beta_{km+l+1})$ and $\HH_{l+1}\HH_l\HH_{l+1}(\beta_{km+l+2})$ by $f\HH_{l+1}\HH_l\HH_{l+1}(\beta_{km+l+1})$ and $\HH_{l+1}\HH_l\HH_{l+1}(\beta_{km+l+2})f^{-1}$ respectively.
\end{Proof}

\subsection{The algebra $\RR$}\label{aR}

In this and the next sections we assume that $t=1$. Let us consider an algebra $\RR$ defined as a path algebra of a quiver with relations. Let $\QQ_{m,n}$ be a quiver whose set of vertices is $\{1,\dots,m\}\times\mathbb{Z}_n$. The set of arrows of $\QQ_{m,n}$ is:
$$\gamma_{i,j}:(i,j)\rightarrow(i+1,j)\mbox{ and }\gamma_{i,j}':(i+1,j)\rightarrow (i,j+1)\,\,(1\le i\le m-1, j\in\mathbb{Z}_n).$$
 Let us define the ideal of relations $\II_{m,n}$. If $m>2$, then it is generated by the elements $\gamma_{i+1,j}\gamma_{i,j}$, $\gamma_{i,j+1}'\gamma_{i+1,j}'$  and $\gamma_{i,j+1}\gamma_{i,j}'-\gamma_{i+1,j}'\gamma_{i+1,j}$ for $1\le i\le m-2$, $j\in\mathbb{Z}_n$. If $m=2$, then the ideal $\II_{m,n}$ is generated by the elements $\gamma_{1,j+1}\gamma_{1,j}'\gamma_{1,j}$ and $\gamma_{1,j+1}'\gamma_{1,j+1}\gamma_{1,j}'$ for $j\in\mathbb{Z}_n$. We define the algebra $\RR$ by the formula $\RR=\kk\QQ_{m,n}/\II_{m,n}$. We omit the indices of $\gamma$ and $\gamma'$ if they can be easily recovered from the context. It is known that the algebra $\RR$ is derived equivalent to the algebra $\NN$.

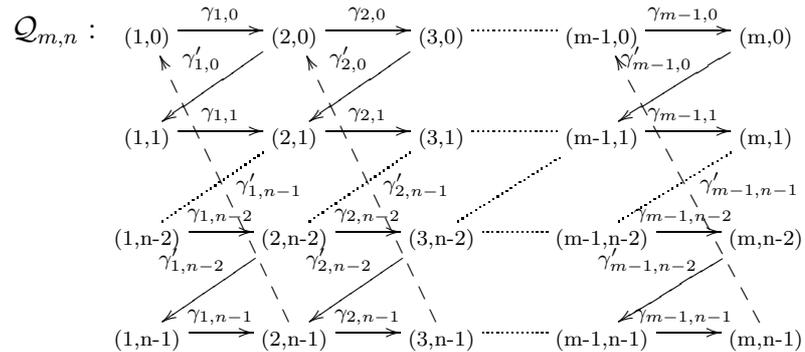
 $$
 \QQ_{m,n}:
\xymatrix
{
\txt{\scriptsize{(1,0)}}\ar@{->}[r]^{\gamma_{1,0}}&\txt{\scriptsize{(2,0)}}\ar@{->}[ld]_{\gamma_{1,0}'}\ar@{->}[r]^{\gamma_{2,0}}&\txt{\scriptsize{(3,0)}}\ar@{->}[ld]_{\gamma_{2,0}'}\ar@{.}[r]&
\txt{\scriptsize{(m-1,0)}}\ar@{->}[r]^{\gamma_{m-1,0}}&\txt{\scriptsize{(m,0)}}\ar@{->}[ld]_{\gamma_{m-1,0}'}\\
\txt{\scriptsize{(1,1)}}\ar@{->}[r]^{\gamma_{1,1}}&\txt{\scriptsize{(2,1)}}\ar@{.}[dl]\ar@{->}[r]^{\gamma_{2,1}}&\txt{\scriptsize{(3,1)}}\ar@{.}[dl]\ar@{.}[r]&\txt{\scriptsize{(m-1,1)}}\ar@{.}[dl]\ar@{->}[r]^{\gamma_{m-1,1}}&\txt{\scriptsize{(m,1)}}\ar@{.}[dl]\\
\txt{\scriptsize{(1,n-2)}}\ar@{->}[r]^{\gamma_{1,n-2}}&\txt{\scriptsize{(2,n-2)}}\ar@{->}[ld]_{\gamma_{1,n-2}'}\ar@{->}[r]^{\gamma_{2,n-2}}&\txt{\scriptsize{(3,n-2)}}\ar@{->}[ld]_{\gamma_{2,n-2}'}\ar@{.}[r]&\txt{\scriptsize{(m-1,n-2)}}\ar@{->}[r]^{\gamma_{m-1,n-2}}&\txt{\scriptsize{(m,n-2)}}\ar@{->}[ld]_{\gamma_{m-1,n-2}'}\\
\txt{\scriptsize{(1,n-1)}}\ar@{->}[r]^{\gamma_{1,n-1}}&\txt{\scriptsize{(2,n-1)}}\ar@{-->}[luuu]_(.45){\gamma_{1,n-1}'}\ar@{->}[r]^{\gamma_{2,n-1}}&\txt{\scriptsize{(3,n-1)}}\ar@{-->}[luuu]_(.45){\gamma_{2,n-1}'}\ar@{.}[r]&\txt{\scriptsize{(m-1,n-1)}}\ar@{->}[r]^{\gamma_{m-1,n-1}}&\txt{\scriptsize{(m,n-1)}}\ar@{-->}[luuu]_(.45){\gamma_{m-1,n-1}'}\\
}
$$

We denote by $e_{i,j}$ the primitive idempotent corresponding to the
vertex $(i,j)$. By $P_{i,j}$ we denote the projective $\RR$-module
$e_{i,j}\RR$ and by $P_{[i_1,j_1][i_2,j_2]}$ the projective
$\RR^{\rm op}\ot\RR$-module $\RR e_{i_1,j_1}\ot e_{i_2,j_2}\RR$.
Also we denote by $\tau$ the automorphism of the algebra $\RR$
defined by the formulas
$$
\tau(e_{i,j})=e_{i,j+1},\,\,\tau(\gamma_{k,j})=\gamma_{k,j+1},\,\,\tau(\gamma_{k,j}')=\gamma_{k,j+1}'\,\,(1\le i\le m,1\le k\le m-1, j\in\mathbb{Z}_n).
$$
Note that $\tau^{-1}$ is a Nakayama automorphism of $\RR$ and so it commutes with any standard equivalence by \cite[Proposition 5.2]{Rickard}.

Let us consider the following $\RR^{\rm op}\ot\RR$-complexes
$$T_i=\left(\bigoplus\limits_{j\in\mathbb{Z}_n}P_{[i,j][i,j]}\xrightarrow{\mu_i} \underline{\RR}\right),$$
where $\mu_i(a\ot b)=ab$ for $a\ot b\in P_{[i,j][i,j]}$. We denote by $\TT_i$ the functor $-\ot_{\RR}T_i$.

The algebra $\RR$ with an additional grading was studied in
\cite{Efimov} (there $\RR$ is
denoted by $A_{m,0}^n$). The arguments from \cite{Efimov} can be applied to our case to obtain the following
result:

\begin{prop}\label{Efimov}
$1)$ $\TT_i$ is an equivalence. In particular, $T_i$ is a two-sided tilting complex for $\RR$.\\
$2)$ There is an injective homomorphism of groups $\lambda:\BB(A_m)\rightarrow \TP(\RR)$ such that $\lambda(s_i)=\TT_i$ ($1\le i\le m$).
\end{prop}

In particular, there is an $A_m$ configuration of $0$-spherical sequences of length $n$ in the category $\Kb\RR.$ For more details on spherical twists see \cite{ST}.

\subsection{Construction of an isomorphism from $\TP(\RR)$ to $\TP(\NN)$}\label{ci}

In this section we  construct explicitly the isomorphism $\TP(\RR)\cong\TP(\NN)$.

Let us consider the following $\RR$-complexes:
$$
U_{i,j}=\left(\underline{P_{1,j}}\xrightarrow{\gamma} P_{2,j}\xrightarrow{\gamma}\dots\xrightarrow{\gamma} P_{i,j}\right).
$$
It can be easily verified that $U:=\oplus_{1\le i\le m,j\in\mathbb{Z}_n}U_{i,j}$ is a tilting complex for $\RR$. For $i\in\mathbb{Z}_{nm}$ we denote by $p(i)$ and $q(i)$ the unique integers such that $0\le p(i)\le m-1$, $0\le q(i)\le n-1$ and $p(i)+mq(i)=i$ in $\mathbb{Z}_{nm}$. Let us define an equivalence $S_{\omega}:\PP_A\rightarrow\add U$  by the formulas
$$
S_{\omega}(P_i)=U_{m-p(i),q(i)},\,\,S_{\omega}(\beta_i)=\begin{cases}
(\underline{\Id_{P_{1,q(i)}}},\dots,\Id_{P_{m-p(i)-1},q(i)}),&\mbox{ if $p(i)\le m-2$,}\\
\gamma_{1,q(i)}'\gamma_{1,q(i)},&\mbox{ if $p(i)=m-1$.}\\
\end{cases}
$$
The tilting functor $(U,\omega)$ defines an isomorphism of groups $$L_{\omega}:\TP(\RR)\rightarrow\TP(\NN)$$ by the formula $L_{\omega}(F)=\bar F_{\omega}FF_{\omega}$, where $\bar F_{\omega}$ is the quasi-inverse equivalence for $F_{\omega}$.
By \cite{Linck} $L_{\omega}(\Pic_0(\RR))=\Pic_0(\NN)$ and $L_{\omega}([1])=[1]$. It is also easy to see that $L_{\omega}(\tau)=\rho^m$.

\begin{Lemma}\label{omega_image} For any projective $\NN$-module $P$ we have
$$
L_{\omega}(\TT_i)(P)\cong\begin{cases}
\HH_{m-i}(P),&\mbox{ if $i>1$,}\\
\QQ_{m-1}(P)[1],&\mbox{ if $i=1$.}\\
\end{cases}
$$
\end{Lemma}
\begin{Proof}
Let us consider the case $i>1$. We have to prove that $F_{\omega}\HH_{m-i}(P_j)\cong \TT_iF_{\omega}(P_j)$ for any $j\in\mathbb{Z}_{nm}$. Let $p(j)$ and $q(j)$ be as above. Note that $2\le p(j)+i\le 2m-1$. Then we have to prove that
$$
F_{\omega}\HH_{m-i}(P_{p(j)+mq(j)})\cong \TT_i(U_{m-p(j),q(j)}).
$$
If $p(j)+i\not\in\{m,m+1\}$, then by definition $F_{\omega}\HH_{m-i}(P_{p(j)+mq(j)})=U_{m-p(j),q(j)}$. If $m-p(j)<i-1$, then it is clear that $\TT_i(U_{m-p(j),q(j)})\cong U_{m-p(j),q(j)}$. If $m-p(j)>i$, then $\TT_i(U_{m-p(j),q(j)})$ is isomorphic to the cone of a map
%$$
%(\underline{\gamma_{i-1,q(j)-1}'},\gamma_{i,q(j)-1}'\gamma_{i,q(j)-1})[2-i]\oplus(\underline{\Id_{P_{i,q(j)}}},\gamma_{i,q(j)})[1-i]
%$$
from
$$
C:=\left(\underline{P_{i,q(j)-1}}\xrightarrow{\Id_{P_{i,q(j)-1}}} P_{i,q(j)-1}\right)[2-i]\oplus \left(\underline{P_{i,q(j)}}\xrightarrow{\Id_{P_{i,q(j)}}} P_{i,q(j)}\right)[1-i]
$$
to $U_{m-p(j),q(j)}$. Since $C\cong 0$, $\TT_i(U_{m-p(j),q(j)})\cong U_{m-p(j),q(j)}$  in $\Kb\RR$. It remains to consider the cases $p(j)=m-i$ and $p(j)=m-i+1$.

Assume that $p(j)=m-i$. In this case $F_{\omega}\HH_{m-i}(P_{p(j)+mq(j)})=F_{\omega}(P_{p(j)+1+mq(j)})=U_{i-1,q(j)}$. At the same time $\TT_i(U_{m-p(j),q(j)})$ is isomorphic to the cone of the map
$$
(\underline{\gamma_{i-1,q(j)-1}'},\gamma_{i,q(j)-1}'\gamma_{i,q(j)-1})[2-i]\oplus \Id_{P_{i,q(j)}}[1-i]
$$
from
$$
C:=\left(\underline{P_{i,q(j)-1}}\xrightarrow{\Id_{P_{i,q(j)-1}}} P_{i,q(j)-1}\right)[2-i]\oplus P_{i,q(j)}[1-i]
$$
to $U_{i,q(j)}$. It is easy to see that this cone is isomorphic to $U_{i-1,q(j)}$.

Let us now consider the case $p(j)=m-i+1$. In this case
$$
F_{\omega}\HH_{m-i}(P_{p(j)+mq(j)})=F_{\omega}(X_{p(j)+mq(j)-1}).
$$
Note that $F_{\omega}(\beta_{p(j)+m(q(j)-1),m})$ is homotopic to
$(\gamma_{i-1,q(j)-1}'\gamma_{i-1,q(j)-1})[2-i]$. Applying $F_{\omega}$ to
$X_{i+mq(j)}$ we obtain the totalization of the following bicomplex
{\footnotesize $$
\xymatrix
{
P_{i,q(j)-1}\ar@/^1pc/@{-->}[rrrrd]^{\gamma'}\\
P_{i-1,q(j)-1}\ar@{->}[rr]^{(-1)^i\Id}\ar@{->}[u]^{\gamma}&&P_{i-1,q(j)-1}\ar@{->}[rr]^{(-1)^{i-1} \gamma'\gamma}&&P_{i-1,q(j)}\\
P_{i-2,q(j)-1}\ar@{->}[rr]^{(-1)^{i-1}\Id}\ar@{->}[u]^{\gamma}&&P_{i-2,q(j)-1}\ar@{->}[u]^{\gamma}&&P_{i-2,q(j)}\ar@{->}[u]^{\gamma}\\
P_{1,q(j)-1}\ar@{->}[rr]^{\Id}\ar@{.}[u]&&P_{1,q(j)-1}\ar@{.}[u]&&\underline{P_{1,q(j)}}\ar@{.}[u]\\
}
$$}
The dotted arrow was constructed according to the algorithm from Section \ref{sde}. It is clear that the totalization is isomorphic to the cone of
$$
P_{i,q(j)-1}[2-i]\xrightarrow{\gamma_{i-1,q(j)-1}'[2-i]} U_{i-1,q(j)}
$$
and it is easy to see that $\TT_i(U_{m-p(j),q(j)})$ is isomorphic to the same cone. The case $i=1$ is analogous and so it is left to the reader.
\end{Proof}

To prove part 2 of Theorem \ref{MainThm} we also need the following lemma.

\begin{Lemma}\label{rot}
$
\QQ_{m-1}^2\HH_{m-2}\dots \HH_{0}(P_i)\cong P_{i-m-1}
$ for all $i\in\mathbb{Z}_{nm}$.
\end{Lemma}
\begin{Proof} If $m\mid i$, then
$$
\QQ_{m-1}^2\HH_{m-2}\dots \HH_{0}(P_i)=\QQ_{m-1}^2\HH_{m-2}\dots \HH_{1}(P_{i+1})=\dots=\QQ_{m-1}^2(P_{i+m-1})=P_{i-m-1}.
$$

Now let $m\mid i-k$ where $1\le k\le m-1$. We have
$$
\begin{aligned}
&\QQ_{m-1}^2\HH_{m-2}\dots \HH_{0}(P_i)=\QQ_{m-1}^2\HH_{m-2}\dots \HH_{k-1}(P_i)
=\QQ_{m-1}^2\HH_{m-2}\dots \HH_{k}(X_{i-1})\\
=&\QQ_{m-1}^2\HH_{m-2}\dots \HH_{k+1}(P_{i-m-1}\xrightarrow{\beta^2} P_{i-m+1}\xrightarrow{\beta^m} \underline{P_{i+1}})\\
=&\dots=\QQ_{m-1}^2(P_{i-m-1}\xrightarrow{\beta^{m-k}} P_{i-k-1}\xrightarrow{\beta^{m}} \underline{P_{i-k-1+m}}).
\end{aligned}
$$
$\QQ_{m-1}(P_{i-m-1}\xrightarrow{\beta^{m-k}} P_{i-k-1}\xrightarrow{\beta^{m}} \underline{P_{i-k-1+m}})$ is the totalization of the bicomplex
$$
\xymatrix
{
P_{i-m-1}\ar@{-->}[rrrrd]^{\beta^{m-k}}\\
P_{i-k-1-m}\ar@{->}[rr]^{\Id_{P_{i-k-1-m}}}\ar@{->}[u]^{\beta^{k}}&&P_{i-k-1-m}\ar@{->}[rr]^{-\beta^{m}}&&\underline{P_{i-k-1}}
}
$$
So
$$
\QQ_{m-1}^2\HH_{m-2}\dots \HH_{0}(P_i)\cong \QQ_{m-1}\left(P_{i-m-1}\xrightarrow{\beta^{m-k}}\underline{P_{i-k-1}}\right)\cong P_{i-m-1}.
$$
\end{Proof}

\subsection{The case $t=1$}\label{case1} In this section we prove part 2 of
Theorem \ref{MainThm}. So we assume that $t=1$.
Since $\RR$ is derived equivalent to $\NN$, it follows from
Corollary \ref{Pic0} that there is some isomorphism
$\zeta_0:\kk^*\rightarrow\Pic_0(\RR)$. First, we prove the following
proposition.

\begin{prop}\label{mono}
There is a monomorphism of groups
$$
\zeta:\AAA_{m,n}\times \kk^*\rightarrow \TP(\RR)
$$
such that:
 $$
 \begin{aligned}
 \zeta(1,a)=\zeta_0(a)\,\,(a\in\kk^*),\,\,\zeta(s_i,1)=\TT_i\,\,(1\le i\le m),\,\,
 \zeta(r_1,1)=[1]\mbox{ and }\zeta(r_2,1)=\tau.
 \end{aligned}
 $$
\end{prop}
\begin{Proof}
To prove the existence of $\zeta$ we have to check that\\
1) every element of $\Pic_0(\RR)$ commutes with $\TT_i$ ($1\le i\le m$), $[1]$ and $\tau$;\\
2) $[1]$ and $\tau$ commute with $\TT_i$ ($1\le i\le m$) and with each other;\\
3) $\TT_i\TT_j=\TT_j\TT_i$ if $1\le i,j\le m$ and $|i-j|>1$;\\
4) $\TT_i\TT_{i+1}\TT_i=\TT_{i+1}\TT_i\TT_{i+1}$ for $1\le i\le m-1$;\\
5) $(\TT_m\dots\TT_1)^{m+1}=\tau^{m+1}[2m]$;\\
6) $\tau^n=\Id_{\RR}$.

Any element of $\Pic_0(\RR)$ is of the form $-\ot_{\RR}
\RR_{\sigma}$ for some $\sigma\in\Aut(\RR)$ such that
$\sigma(e_{i,j})=e_{i,j}$ for all $1\le i\le m$, $j\in\mathbb{Z}_n$
(here $\RR_{\sigma}$ is an $\RR$-bimodule which coincides with $\RR$
as a left module and has the right multiplication $*$ defined by the
formula $a*b=a\sigma(b)$ for $a\in \RR_{\sigma}$, $b\in\RR$). The
map
$\psi:\RR\ot_{\RR}\RR_{\sigma}\rightarrow\RR_{\sigma}\ot_{\RR}\RR$
defined by the formula $\psi(a\ot b)=a\ot\sigma^{-1}(b)$ is an
isomorphism of bimodules. Also for all $1\le i\le m$ and
$j\in\mathbb{Z}_n$ we have an isomorphism of bimodules
$\psi_{i,j}:P_{[i,j],[i,j]}\ot_{\RR}\RR_{\sigma}\rightarrow
\RR_{\sigma}\ot_{\RR}P_{[i,j],[i,j]}$ defined by the formula
$\psi_{i,j}(a\ot b\ot c)=a\ot e_{i,j}\ot \sigma^{-1}(bc)$. Then the
map
$\left(\oplus_{j\in\mathbb{Z}_n}\psi_{i,j},\underline{\psi}\right):T_i\ot_{\RR}
\RR_{\sigma}\rightarrow\RR_{\sigma}\ot_{\RR} T_i$ is an isomorphism
of complexes of $\RR$-bimodules. It is clear that the shift commutes
with any standard equivalence. As was mentioned above  $\tau$
commutes with any standard equivalence too. So 1) and 2) hold. 3)
and 4) follow from Proposition \ref{Efimov} and 6) is clear.

By \cite[Lemma 4.30]{Efimov} (applied to our case) we have
\begin{equation}\label{center_eq}
(\TT_m\dots\TT_1)^{m+1}(P_{i,j})=P_{i,j+m+1}[2m].
\end{equation}
So the element $(\TT_m\dots\TT_1)^{m+1}[-2m]\tau^{-(m+1)}$ lies in
$\Pic_0(\RR)$. There is a canonical map from $\TP(\RR)$ to the group
of isomorphism classes of $\RR$-bimodules inducing stable
autoequivalences of $\RR$ (see \cite[Section 3.4]{RouZim} for
details). This map is injective on $\Pic_0(\RR)$. At the same time
it sends $\TT_i$ to $\RR$ and $[-2m]$ to $\Omega^{2m}_{\RR^{\rm
op}\ot\RR}(\RR)$. So it is enough to prove that
$\Omega^{2m}_{\RR^{\rm op}\ot\RR}(\RR)$ is isomorphic to
$\RR_{\tau^{m+1}}$. From \cite[Section 4]{ErdHolm} it follows that
$\Omega^{2m}_{\NN^{\rm op}\ot\NN}(\NN)\cong \NN_{\rho^{m(m+1)}}$.
There are $\NN$--$\RR$-bimodule $M$ and $\RR$--$\NN$-bimodule $\bar
M$ which correspond to $F_{\omega}$ and $\bar F_{\omega}$
respectively (see \cite[Corollary 2.15]{RouZim}). Then we have
$$\Omega^{2m}_{\RR^{\rm op}\ot\RR}(\RR)\cong \bar M\ot_{\NN}\Omega^{2m}_{\NN^{\rm op}\ot\NN}(\NN)\ot_{\NN}M\cong \bar M\ot_{\NN}\NN_{\rho^{m(m+1)}}\ot_{\NN}M\cong \RR_{\tau^{m+1}}$$
in the stable category of $\RR$-bimodules. Here we use the fact that $\rho^{-m}$ is a Nakayama automorphism and so $\rho^{m(m+1)}F_{\omega}\cong F_{\omega}\tau^{m+1}$  by \cite[Proposition 5.2]{Rickard}.

Let us now consider an element
$(sr_1^{-l_1}r_2^{-l_2},a)\in\AAA_{m,n}\times \kk^*$, where
$s\in\BB(A_n)$. Suppose that
$\zeta(sr_1^{-l_1}r_2^{-l_2},a)=\Id_{\Kb \RR}$. Then $\zeta(s)$
commutes with $\TT_i=\zeta(s_i)$ for any $1\le i\le m$. Since
$\zeta|_{\BB(A_n)}$ is injective by Proposition \ref{Efimov}, $s$
belongs to the center of $\BB(A_n)$. Then by results of
\cite{BerSait} we have $s=(s_n\dots s_1)^{k(m+1)}$ for some
$k\in\mathbb{Z}$. So
$(\TT_m\dots\TT_1)^{k(m+1)}=\zeta_0(a)\tau^{l_2}[l_1]$. Then it
follows from \eqref{center_eq} that $l_1=2mk$ and $n\mid
l_2-k(m+1)$. Thus, $sr_1^{-l_1}r_2^{-l_2}=1$, $\zeta_0(a)=1$. Since
$\zeta_0:\kk^*\rightarrow\Pic_0(\RR)$ is an isomorphism, $a=1$ and
so $(sr_1^{-l_1}r_2^{-l_2},a)=(1,1)$. Consequently, $\zeta$ is a
monomorphism.
\end{Proof}

\begin{Proof}[Proof of part 2 of Theorem \ref{MainThm}] By Proposition \ref{mono} it is enough to prove that $L_{\omega}(\TT_i)$ ($1\le i\le m$), $L_{\omega}(\Pic_0(\RR))$, $L_{\omega}(\tau)$ and $[1]$ generate $\TP(\NN)$.
Denote by $X$ the subgroup of $\TP(\NN)$ generated by these
elements. $\Pic_0(\NN)$ belongs to $X$. By Lemma \ref{omega_image} $\HH_i$ ($0\le i\le m-2$),
$\QQ_{m-1}$ belong to $X$. By Lemma \ref{rot} $\rho^{-(m+1)}$ belongs to $X$, thus by Proposition
\ref{Picard} and since $\rho^{m} \in X$, the subgroup $\Pic(\NN)$
also belongs to $X$. Since
$$
\HH_{m-1}=\rho\HH_{m-2}\rho^{-1}\mbox{ and
}\QQ_i=\rho^{i+1-m}\QQ_{m-1}\rho^{m-1-i}\,\,(0\le i\le m-2),
$$
$\HH_{m-1}$ and $\QQ_i$ ($0\le i\le m-2$) belong to $X$. Thus by
part 3 of Theorem \ref{Nakayama_gen} we have $X=\TP(\NN)$.
\end{Proof}

From here on we assume that $t\ge 2$. Let us denote $\NN(ntm,tm)$ by $\td\NN$. Moreover we denote by
$\td\HH_l,\td\QQ_l:\Kb\td\NN\rightarrow\Kb\td\NN$
($l\in\mathbb{Z}_{tm}$) the standard equivalences whose construction
for $\td\NN$ coincides with that of $\HH_l$ and $\QQ_l$
($l\in\mathbb{Z}_m$) for $\NN$ respectively. Similarly, we define
$\td e_i$, $\td\beta_i$, $\td P_i$ ($i\in\mathbb{Z}_{ntm}$), $\td
X_i$, $\td Y_i$ ($i\in\mathbb{Z}_{tm}$).

\subsection{Some properties of $\BB(\td A_{N-1})$}\label{BAff}

In this section we list some properties of the group $\BB(\td
A_{N-1})$. Most of them are known but we collect them here for
convenience. Let us introduce the notation $\Pi:=s_N\dots s_1$. The
following result can be found in \cite{CharPei}:

\begin{Lemma}\label{affine_prop} The group $\BB(\td A_{N-1})$ is torsion-free and its center is trivial.
\end{Lemma}

The next fact can be found in \cite[Section 5.1]{Intan}:

\begin{Lemma}\label{td_A_A_mono}
There is a monomorphism of groups $\phi_N:\BB(\td A_{N-1})\rightarrow \BB(A_N)$ such that
$$\phi_N(s_i)=s_{i+1}\,\,(0\le i\le N-2)\mbox{ and }\phi_N(s_{N-1})=s_N\Pi s_1 \Pi^{-1}s_N^{-1}.$$
\end{Lemma}

\begin{coro}\label{faith}
There is a monomorphism of groups $\td\eta:\BB(\td
A_{tm-1})\rightarrow \TP(\td\NN)$ such that $\td\eta(s_i)=\td\HH_i$
for $i\in\mathbb{Z}_{tm}$.
\end{coro}
\begin{Proof}
By Proposition \ref{mono} there is a monomorphism of groups
$\eta':\BB(A_{tm})\rightarrow \TP(\td\NN)$ which sends $s_i$ to
$L_{\omega}(\TT_{tm-i+1})$. Let us compute $\eta'\phi_{tm}(s_i)$. By
Lemma \ref{omega_image} we have
$$
\begin{aligned}
\td\HH_i&=L_{\omega}(\TT_{tm-i})W_i=\eta'(s_{i+1})W_i\,\,(0\le i\le tm-2),\\
\td\QQ_{tm-1}&=L_{\omega}(\TT_1)[-1]W_{tm-1}=\eta'(s_{tm})[-1]W_{tm-1}
\end{aligned}
$$
for some $W_i\in\Pic_0(\td\NN)$, for $i\in\mathbb{Z}_{tm}$. By Lemma
\ref{rot} we have
$$
\td\QQ_{tm-1}^2\td\HH_{tm-2}\dots \td\HH_0=\rho^{-(tm+1)}W
$$
for some $W\in\Pic_0(\td\NN)$.  Since $\Pic_0(\td\NN)$ lies in the
center of $\TP(\td\NN)$, we have
$$
\eta'(s_{tm}\Pi s_1\Pi^{-1}s_{tm}^{-1})W_0=\rho^{-(tm+1)}\td\HH_0\rho^{tm+1}=\td\HH_{tm-1}.
$$
By Lemma \ref{aff_braid} there is a homomorphism $\td\eta$ with the
required properties. If $\td\eta(s)=\Id_{\Kb\td\NN}$ for some
$s\in\BB(\td A_{tm-1})$, then it is easy to see that
$\eta'\phi_{tm}(s)\in\Pic_0(\td\NN)$. Since
$\Pic_0(\td\NN)\cap\Im\eta'=\{\Id_{\Kb\td\NN}\}$ (this follows, for
example, from Proposition \ref{mono}) and the maps $\eta'$ and
$\phi_{tm}$ are monomorphic by Lemma \ref{td_A_A_mono}, we have
$s=1$. Consequently, $\td\eta$ is a monomorphism.

\end{Proof}

\begin{Lemma}\label{cylmono}
There is a monomorphism of groups $\psi:\BB(\td A_{m-1})\rightarrow \BB(\td A_{tm-1})$ defined by the formula
$$\psi(s_i)=s_is_{i+m}\dots s_{i+(t-1)m}\,\,(i\in\mathbb{Z}_m).$$
\end{Lemma}
\begin{Proof}It is easy to check that $\psi$ is a homomorphism of groups. We only
need to check that it is monomorphic.

We are going to consider $\BB(\td A_{m-1})$ as a subgroup of the braid group of a two-dimensional surface. Let
$L$ be the complex plane without the origin $\mathbb{C}^*$, or
equivalently a sphere with two punctures, denote by
$\Delta_m=\{(x_0, \dots, x_{m-1}) \in L^m \mid x_i=x_j \mbox{ for
some } i \neq j\}.$ The symmetric group $S_m$ acts on $L^m$, the
action is free on $L^m \setminus \Delta_m$. Consider $X_m=(L^m
\setminus \Delta_m)/S_m$, fix a base point $(b_0, \dots, b_{m-1})
\in X_m$, then the $m$-strand braid group of $L$ is defined to be
the fundamental group of $X_m$, we will denote it by $Br_m(L)$.
Each element of $Br_m(L)$ may be represented by a set of paths
$\gamma_0, \dots, \gamma_{m-1}: [0,1] \rightarrow L$ such that there
is a permutation $\pi$ of $\{0, \dots, m-1\}$ such that
$\gamma_i(0)=b_i, \gamma_i(1)=b_{\pi(i)}$ and for each $t \in [0,1]$
the points $\gamma_i(t)$ are all distinct. Then $\BB(\td A_{m-1})$
is a normal subgroup in $Br_m(L)$ with the quotient isomorphic to
$\mathbb{Z}$. A good choice for a set of basepoints for $m$-strand
braids in $\mathbb{C}^*$ is the set of $m$-th roots of unity. Then
the standard generators $s_i$ correspond to the standard braid
generators $s'_i$, which exchange two adjacent basepoints in the
simplest possible way \cite{Allcock}.

Consider a $t$-sheeted covering of $L$ with the covering space
$\tilde{L}$ homeomorphic to $L$, $p: \tilde{L} \rightarrow L$. If we
identify $L$ with $\mathbb{C}^*$ this covering can be given by
$z^t$. Consider also the group $Br_{tm}(\tilde{L})$. The base point
$(\tilde{b}_0, \dots, \tilde{b}_{tm-1}) \in \tilde{X}_{tm}$ is
chosen in such a way that $p(\tilde{b}_{i+sm}) = b_i$, $s=0, \dots,
t-1$. This covering induces a covering $\tilde{L} \setminus
\{\tilde{b}_0, \dots, \tilde{b}_{tm-1}\} \rightarrow L \setminus
\{b_0, \dots, b_{m-1}\},$ which will be also denoted by $p.$ There is
a lift of $s'_i$ from $L$ to $\tilde{L}$, it will be denoted by
$\tilde{s}'_i$. The element $\tilde{s}'_i$ exchanges $\tilde{b}_{i+sm}$ and
$\tilde{b}_{i+1+sm}$ for $s=0,\dots,t-1.$ The homomorphism $\psi$ is consistent with the
lift $\tilde{s}'_i$ of $s'_i$.

The fundamental group $\pi_1(\tilde{L} \setminus \{\tilde{b}_0,
\dots, \tilde{b}_{tm-1}\})$ is isomorphic to $F_{tm+1}$. Since $p$ is
a covering it induces a monomorphism $$p_* : \pi_1(\tilde{L}
\setminus \{\tilde{b}_0, \dots, \tilde{b}_{tm-1}\}) \simeq F_{tm+1}
\rightarrow \pi_1(L \setminus \{b_0, \dots, b_{m-1}\}) \simeq
F_{m+1}.$$ The index of $F_{tm+1}$ in $F_{m+1}$ is equal to $t.$

Let us check that $\BB(\td A_{m-1})$ acts faithfully on $\pi_1(L
\setminus \{b_0, \dots, b_{m-1}\}) \simeq F_{m+1}$. Let us consider the
mapping class group $MCG_m(L)$ of $L \setminus \{b_0, \dots,
b_{m-1}\}$, i.e. the group of all path components in the space of all
homeomorphisms from $L$ to $L$ which keep $\{b_0, \dots, b_{m-1}\}$
fixed as a set. By \cite{Birman} there is a long exact sequence
$$\dots \rightarrow \pi_1(Homeo(L)) \rightarrow Br_m(L)
\xrightarrow{w} MCG_m(L) \rightarrow \pi_0(Homeo(L)) \rightarrow
1,$$ where $Homeo(L)$ is the space of all homeomorphisms from $L$ to
$L$. The group $\pi_1(Homeo(L))$ is isomorphic to $\mathbb{Z}$, its
image in $Br_m(L)$ is generated by the central element of $Br_m(L)$,
corresponding to the rotation of $\mathbb{C}^*$ by $2\pi$. Thus the
intersection of $\BB(\td A_{m-1})$ and the kernel of $w$ is trivial
and $w$ induces an embedding from $\BB(\td A_{m-1})$ to $MCG_m(L)$.
One can consider the stabilizer of one puncture in $MCG_m(L)$, if we
identify $L$ with $\mathbb{C}^*$ it would be the stabilizer of the
origin, this stabilizer is a subgroup of the mapping class group
$MCG_{m+1}(\mathbb{C})$ of all homeomorphisms $\mathbb{C}
\rightarrow \mathbb{C}$ which keep $\{b_0, \dots, b_{m-1}, 0\}$ fixed as
a set. It is easy to see that $\BB(\td A_{m-1})$ stabilizes the
puncture, thus it can be embedded into $MCG_{m+1}(\mathbb{C})$. It
is well known that $MCG_{m+1}(\mathbb{C})$ coincides with the
$m+1$-strand braid group of $\mathbb{C}$ or equivalently a sphere
with one puncture and that this group acts faithfully on $\pi_1(L
\setminus \{b_0, \dots, b_{m-1}\}) \simeq F_{m+1}$. Thus $\BB(\td
A_{m-1})$ is embedded into a group which acts faithfully on $\pi_1(L
\setminus \{b_0, \dots, b_{m-1}\}) \simeq F_{m+1}$, hence it acts
faithfully on $\pi_1(L \setminus \{b_0, \dots, b_{m-1}\}) \simeq
F_{m+1}$.

If we denote by $y_i$ the loop around $b_i$ and by $y$ the loop around the origin, then the action of  $\BB(\td A_{m-1})$ is defined as follows: $$s_i(y_i)=y_{i+1},\,s_i(y_{i+1})=y_{i+1}^{-1}y_iy_{i+1},\,s_i(y_j)=y_j\,(j \neq i, i+1),\,s_i(y)=y.$$ Similarly, $\BB(\td A_{tm-1})$ acts faithfully on $\pi_1(\tilde{L} \setminus
\{\tilde{b}_0, \dots, \tilde{b}_{tm-1}\}) \simeq F_{tm+1}$ and these two actions are compatible, in particular,
$p_*\psi(g)(\gamma)$ is homotopic to $g(p_*(\gamma))$, for $g \in
\BB(\td A_{m-1})$, ${\gamma\in\pi_1(\tilde{L} \setminus
\{\tilde{b}_0, \dots, \tilde{b}_{tm-1}\})}$, this can be easily checked on the basis.

Assume that $\psi$ is not monomorphic, i.e. there exists $g \neq 1
\in \BB(\td A_{m-1})$ such that $\psi(g)=1$. Since $g \neq 1$ and
the action of $\BB(\td A_{m-1})$ on $F_{m+1}$ is faithful $g(\gamma) \neq
\gamma$ for some $\gamma \in F_{m+1}$, but $\gamma^t \in
p_*(F_{tm+1})$, hence $\gamma^t=p_*(\gamma')$ for some $\gamma' \in
F_{tm+1}$ and $$g(\gamma)^t = g(\gamma^t)=g(p_*(\gamma'))=
p_*(\psi(g)(\gamma'))=p_*(\gamma')=\gamma^t.$$ Finally, $g(\gamma)
\neq \gamma$ and $g(\gamma)^t = \gamma^t,$ which is not possible in
a free group \cite[Chapter I, Propositions 2.16, 2.17]{LynSch}. \end{Proof}

It follows from Corollary \ref{faith} and Lemma \ref{cylmono} that
there is a monomorphism of groups $\bar\eta=\td\eta\psi:\BB(\td
A_{m-1})\rightarrow \TP(\td\NN)$ which sends $s_i\in\BB(\td
A_{m-1})$ ($0\le i\le m-1$) to $\HH_i\HH_{i+m}\dots \HH_{i+(t-1)m}$.

\subsection{Faithful action of $\BB(\td A_{m-1})$}\label{fa}

Let us consider the homomorphism $\Delta:C_t\rightarrow
\Aut(\td\NN)$ defined by the formula $\Delta(r)=\td\rho^{nm}$ (where
$\td\rho$ is the rotation of $\td\NN$).

For $l\in\mathbb{Z}_m$ let us introduce the $\td\NN$-complexes
$$
\bar H_l=\Big(\bigoplus_{i\in\mathbb{Z}_{ntm}, m\nmid i-l}\td P_i\Big)\oplus\Big(\bigoplus_{i\in\mathbb{Z}_{ntm}, m\mid i-l}\td X_i\Big).
$$
In this case we define an equivalence $S_{\bar\theta_l}:\PP_{\td\NN}\rightarrow \add \bar H_l$ in the following way:
$$
S_{\bar\theta_l}(\td P_i)=\begin{cases}
\td P_i,&\mbox{ if $m\nmid i-l$ and $m\nmid i-1-l$,}\\
\td P_{i+1},&\mbox{ if $m\mid i-l$,}\\
\td X_{i-1},&\mbox{ if $m\mid i-1-l$,}
\end{cases}
$$
$$
S_{\bar\theta_l}(\td\beta_i)=\begin{cases}
\td\beta_i,&\mbox{ if $m\nmid i-l+1$ and $m\nmid i-1$,}\\
\td\beta_{i+1}\td\beta_i,&\mbox{ if $m\mid i-l+1$,}\\
\Id_{\td P_{i+1}},&\mbox{ if $m\mid i-l$.}
\end{cases}
$$
It is clear that $F_{\bar\theta_l}=\td\HH_l\td\HH_{l+m}\dots
\td\HH_{l+(t-1)m}$. Note that for any $i\in\mathbb{Z}_{ntm}$ there
is an isomorphism $f_i:\td P_i\rightarrow \td P_{i+nm}\#r$ defined
by the equality $f_i(\td e_i)=\td e_{i+nm}\#r$. We define
$s_{i,l}:S_{\bar\theta_l}(\td P_i)\rightarrow S_{\bar\theta_l}(\td
P_{i+nm})\#r$ in the following way. If $S_{\bar\theta_l}(\td
P_i)=\td P_j$ for some $j$, then $s_{i,l}=f_j$. If
$S_{\bar\theta_l}(\td P_i)=\td X_i$, then
$s_{i,l}=(f_{i-tm-1},f_{i-tm},\underline{f_{i}})$. Let
$\pi_{i,l}:\bar H_l\rightarrow S_{\bar\theta_l}(\td P_i)$ and
$\iota_{i,l}:S_{\bar\theta_l}(\td P_i)\rightarrow \bar H_l$
($i\in\mathbb{Z}_{ntm}$) be the canonical projection and inclusion.
We define $\psi_l\in\Kb\NN(\bar H_l,\bar H_l\#r)$ by the equality
$$
\psi_l=\sum_{i\in\mathbb{Z}_{ntm}}(\iota_{i+nm,l}\#r)s_{i,l}\pi_{i,l}
$$
Then $(\bar H_l,\bar\theta_l,\psi_l)$ is a tilting $C_t$-functor from $\td\NN$ to itself.

$\td\NN\#C_t$ is Morita equivalent to $\NN$.
This equivalence can be constructed in such a way that the induced
isomorphism $L:\TP(\NN)\cong\TP(\td\NN\#C_t)$ satisfies the
condition $L(\HH_i)=F_{\bar\theta_i\#\psi_i}$ ($i\in\mathbb{Z}_m$)
(see \cite[Lemma 3]{VolkZvo}). It is also clear that
$L(\Pic(\NN))=\Pic(\td\NN\#C_t)$.

\begin{prop}\label{affine_faith}
The homomorphism $\eta:\BB(\td A_{m-1})\rightarrow \TP(\NN)$ defined in Lemma \ref{aff_braid} is injective.
\end{prop}
\begin{Proof}
Recall that there is a homomorphism $\bar\eta=\td \eta \psi:\BB(\td A_{m-1}) \rightarrow \TP(\td\NN)$
such that $\bar\eta(s_i) = F_{\bar \theta_i}$. The homomorphism $\psi$ defined in Lemma \ref{cylmono}
is injective. The homomorphism $\td \eta$ is also injective by Corollary \ref{faith}. Hence
$\bar\eta$ is injective. By Lemma \ref{group_act} there are such
$W_i\in\TP_{C_t}(\td\NN)$ ($i\in\mathbb{Z}_m$) that
$\Phi_{\td\NN}(W_i)=F_{\bar\theta_i}$ and
$\Psi_{\td\NN}(W_i)=L(\HH_i)$, i.e.
$\Phi_{\td\NN}(W_i)=\bar\eta(s_i)$ and
$\Psi_{\td\NN}(W_i)=L\eta(s_i)$. Suppose that $s=s_{i_1}\dots
s_{i_k}\in\BB(\td A_{m-1})$ is such that $\eta(s)=\Id_{\Kb\NN}$.
Then $\Psi_{\td\NN}(W_{i_1}\dots W_{i_k})=\Id_{\Kb(\td\NN\#C_t)}$
and so
$$\bar\eta(s)=\Phi_{\td\NN}(W_{i_1}\dots W_{i_k})\in\Pic(\td\NN)$$
by Lemma \ref{group_act}. By Proposition \ref{Picard} we have
$\bar\eta(s)=\td\rho^la$ for some integer number $l$ and
$a\in\Pic_0(\td\NN)$. Then $\bar\eta(s^{ntm})=a^{ntm}$ lies in the
center of $\TP(\td\NN)$. Since $a^{ntm}\in\bar\eta(\BB(\td
A_{m-1}))\cong\BB(\td A_{m-1})$, we have
$\bar\eta(s^{ntm})=\Id_{\Kb\td\NN}$ by Lemma \ref{affine_prop}. Then
$s^{ntm}=1$ and by the same lemma we have $s=1$.
\end{Proof}

\subsection{The case $t>1$}\label{case>1}

In this section we prove the remaining part of Theorem
\ref{MainThm}.

\begin{Proof}[Proof of part 3 of Theorem \ref{MainThm}]
Let us denote $\floor{\frac{t+n-1}{n}}$ by $N$. Let us define a
homomorphism $\zeta:\big(\BB(\td A_{m-1})\rtimes_{\varphi_{m,n}}
C_{nm}\big)\times\SSS_N(\kk)\times C_{\infty}\rightarrow\TP(\NN)$ on
the generators:
$$
\begin{aligned}
\zeta(s_i,1,1,1)&=\HH_i\,\,(i\in \mathbb{Z}_m),\,\,\zeta(1,r,1,1)=\rho,\,\,\zeta(1,1,1,r)=[1]\\
\zeta(1,1,c,1)&=\mu_c\,\,(c\in\SSS_N(\kk)).
\end{aligned}
$$
Let us prove that $\zeta$ is well defined. Propositions \ref{Picard}
and \ref{affine_faith} together with the fact that the shift
commutes with any standard equivalence yield that it is enough to
prove that
\begin{enumerate}
\item $\rho \HH_i\cong \HH_{i+1}\rho$ ($i\in\mathbb{Z}_m$);
\item $\mu_c\HH_i\cong \HH_i\mu_c$ ($i\in\mathbb{Z}_m$, $c\in\SSS_N(\kk)$).
\end{enumerate}
The first equality is obvious. Let us prove the second equality.
Since $\mu_c\HH_i=\rho^{-i}\mu_c\HH_0\rho^i$
 and
$\HH_i\mu_c=\rho^{-i}\HH_0\mu_c\rho^i,$
 it is enough to prove the second equality for $i=0$.

The equivalence $\HH_0\mu_c$ corresponds to the equivalence
$S:\PP_{\NN}\rightarrow\add H_0$ which coincides with $S_{\theta_0}$
on $P_i$ for $i\in\mathbb{Z}_{nm}$ and on
$\beta_i$ for $i\in\mathbb{Z}_{nm}\setminus\{0\}$ and is defined on
$\beta_0$ by the formula
$S(\beta_0)=u_c=S_{\theta_0}(\beta_0)u_c:P_1\rightarrow X_0$. On the
other hand the equivalence $\mu_c \HH_0$ corresponds to the
equivalence $\td S:\PP_{\NN}\rightarrow\add H_0$ which coincides
with $S_{\theta_0}$ on $P_i$ for $i\in\mathbb{Z}_{nm}$ and on
$\beta_i$ for $i\in\mathbb{Z}_{nm}\setminus\{-1\}$ and is defined on
$\beta_{-1}$ by the formula $\td
S(\beta_{-1})=u_cS_{\theta_0}(\beta_{-1}):P_{-1}\rightarrow P_1$.
Then the collection of morphisms $\alpha_i:S(P_i)\rightarrow \td
S(P_i)$ defined by the equalities $\alpha_i=\Id_{S(P_i)}$
($i\in\mathbb{Z}_{nm}\setminus\{0\}$) and
$\alpha_0=u_c:P_1\rightarrow P_1$ defines an isomorphism of functors
$\alpha:S\rightarrow \td S$. So we prove that the definition of
$\zeta$ is correct.

It follows from part 2 of Theorem \ref{Nakayama_gen} and Proposition
\ref{Picard} that $\zeta$ is surjective. So it remains to prove that
it is injective. Assume that
$\zeta(s,r^{-k_1},c,r^{-k_2})=\Id_{\Kb\NN}$. Then
\begin{equation}\label{eq_1}
\zeta(s^{nm},1,1,1)=\zeta(1,1,c^{-nm},r^{k_2nm}).
\end{equation}
Let us denote $\zeta(1,1,c^{-nm},r^{k_2nm})$ by $W$. By \eqref{eq_1}
the element W belongs to $G=\zeta\big(\BB(\td A_{m-1})\big)$. On the
other hand $W$ commutes with any element of $G$ and $G\cong \BB(\td
A_{m-1})$ by Proposition \ref{affine_faith}. Then
$\zeta(s,1,1,1)^{nm}=W=\Id_{\Kb\NN}$ by Lemma \ref{affine_prop}. By
the same lemma $\zeta(s,1,1,1)=\Id_{\Kb\NN}$ and so $s=1$ by
Proposition \ref{affine_faith}. It is clear that
$\zeta(1,r^{-k_1},c,r^{-k_2})=\Id_{\Kb\NN}$ iff $nm\mid k_1$,
$c=1_{\SSS_N(\kk)}$ and $k_2=0$. Thus, the injectivity of
$\zeta$ is proved.
\end{Proof}

\end{document}